
\input amstex
\documentstyle{amsppt}
\topmatter
\title
The 2-dimensional jacobian conjecture via Klein's program
\endtitle
\author P.Katsylo
\endauthor
\address
Independent University of Moscow
\endaddress
\email katsylo\@katsylo.msk.ru
\endemail
\date December 1, 1995
\enddate
\abstract
We investigate the 2-dimensional jacobian conjecture
via Klein's program.
\endabstract
\endtopmatter
\document
\magnification 1200
\define\C{\Bbb C}
\define\A{\Bbb A}
\bf \S0. \rm
F. Klein noticed that we can reformulate some algebraic
geometrical problems into invariant theoretical problems
and then try to solve these invariant theoretical problems.
Note that Klein's program was considered as pure speculative idea
(no one algebraic geometrical problem was solved by
Klein's program). The purpose of our investigation is to
develop Klein's program for algebraic geometrical
problems. The 2-dimensional Jacobian conjecture is taken as
test problem for the development of Klein's program.
This preliminary paper contains a discription of two
approaches to the 2-dimensional Jacobian conjecture via
Klein's program.
\par
Klein's program consists of the following two steps.
\roster
\item A reformulation of an algebraic geometrical
problem into an invariant theoretical problem.
\item A solution of the invariant theoretical problem.
\endroster
Note that usually there are many invariant theoretical
reformulations of an algebraic geometrical problem.
\par
Let us recall the Jacobian conjecture (see \cite{2}, \cite{3}).
\proclaim{The Jacobian conjecture} Suppose the jacobian of a polynomial
mapping
$$F : \C^n \rightarrow \C^n
$$
is equal to 1; then $F$ is an isomorphism.
\endproclaim
In \S1 we formulate two general invariant theoretical
problems. Some invariant theoretical reformulations of
some algebraic geometrical problems are partial cases
of these invariant theoretical problems. In \S2 we define
covariant $Q$ (quasiresultant). In \S3 we recall some facts
about representations of the groups $SL_2$ and $SL_3$.
In \S4 we define group $G$, some $G$-moduli, and covariants.
In \S5 we formulate conjecture 1. This conjecture is true
iff the 2-dimensional Jacobian conjecture is true. Conjecture 1
is a partial case of general invariant theoretical problem *
from \S1. In \S6 we formulate conjecture 2. This conjecture is true
iff the 2-dimensional Jacobian conjecture is true. Conjecture 1
is a partial case of general invariant theoretical problem **
from \S1.
\linebreak
\bf Remark 0.1. \rm
Let us describe
Klein's program for the 2-dimensional
Jacobian conjecture.
Let $V$ be the linear space of pairs
$(F_1,F_2)$ of polynomials of degree $\le n$
in the variables $z_1, z_2$,
$$
X_{\neq 0} = \{(F_1,F_2) \in V |\ the\ jacobian\ of\
(F_1,F_2)\ is\ equal\ to\ a\ nonzero\ constant\},
$$
$$
X_0 = \{(F_1,F_2) \in V |\ the\ jacobian\ of\
(F_1,F_2)\ is\ equal\ to\ 0\},
$$
$I(X_{\neq 0}) \subset \C[V]$ be the ideal
of $X_{\neq 0}$, and
$I(X_0) \subset \C[V]$ be the ideal of $X_0$.
There exist elements $s \in \C[V]$
such that
$s \in I(X_{\neq 0}) \setminus I(X_0)$ iff
the 2-dimensional Jacobian conjecture
is true for polynomial mappings of
polynomial degree $\le n$.
At the present moment there exists only one tool
(Puiseux expansions) to construct
elements of $I(X_{\neq 0}) \setminus I(X_0)$
(these elements give, for example,
Folk-theorem mentioned in \S6).
Klein's program is the way to obtain
elements of $I(X_{\neq 0}) \setminus I(X_0)$
by means of invariant theory.
(Even simple invariant theoretical
constructions give some elements of
$I(X_{\neq 0}) \setminus I(X_0)$, but these elements
do not give a solution of the 2-dimensional
Jacobian conjecture.)
\par
We use in the paper some facts of representation theory
of semisimple linear algebraic groups
(see \cite{1}).
\bigpagebreak
\par
\bf \S1. \rm
In this section we formulate two general invariant
theoretical problems.
\par
For polynomial mappings
$\alpha , \beta : V \rightarrow W$
we write $\left.\alpha(v)\right|_{\beta(v)=0} = 0$
iff $v \in V, \beta(v)=0$ imply
$\alpha(v)=0$.
\proclaim{Problem *} Let $G$ be a reductive linear algebraic
group, $G : V, U, \widetilde U$ be linear representations,
$$
\varphi : V \rightarrow U, \quad
\widetilde {\varphi} : V \rightarrow \widetilde U
$$
be homogeneous covariants.
Is it true that
$\left.\varphi(v)\right|_{\widetilde {\varphi} (v) = 0} = 0$?
\endproclaim
\proclaim{Definition 1.1} Let $G$ be a linear algebraic
group, $G : V, U, \widetilde U$ be linear representations,
$$
\varphi : V \rightarrow U, \quad
\widetilde {\varphi} : V \rightarrow \widetilde U
$$
be covariants. A covariant
$$
\gamma : V \times U \rightarrow W
$$
is called $(\varphi, \widetilde \varphi)$-identity if
$$
\left.\gamma(v,\varphi(v)) \right|_{\widetilde {\varphi}(v) = 0} = 0.
$$
\endproclaim
\proclaim{Problem **} Let $G$ be a reductive linear algebraic
group, $G : V, U, \widetilde U$ be linear representations,
$$
\varphi : V \rightarrow U, \quad
\widetilde {\varphi} : V \rightarrow \widetilde U
$$
be homogeneous covariants, and
$v_0 \in V, u_0 \in U$.
Does exist a bihomogeneous
$(\varphi, \widetilde \varphi)$-identity $\gamma$ such
that $\gamma(v_0,u_0) \neq 0$?
\endproclaim
\bigpagebreak
\par
\bf \S2. \rm
In this section we define covariant $Q$.
\par
Consider the group $SL_m$. The group $SL_m$ acts canonically
in the spaces $\C^m$, $\C^{m\ast}$, $S^n\C^{m\ast}$, \dots.
\par
Consider the resultant
$$
R : S^{n_1}\C^{m\ast} \times \dots \times  S^{n_m}\C^{m\ast}
\rightarrow \C.
$$
Set $N = n_1 \cdot \dots \cdot n_m$. The resultant $R$ is a
polyhomogeneous (of polydegree
$(\frac{N}{n_1}$, \dots, $\frac{N}{n_m})$) \ $SL_m$-covariant.
Consider the case $n_m = 1$. In this case the resultant $R$
defines canonically the polyhomogeneous
(of polydegree $(\frac{N}{n_1}$, \dots, $\frac{N}{n_{m-1}})$) \quad
$SL_m$-covariant
$$
Q :  S^{n_1}\C^{m\ast} \times \dots \times  S^{n_{m-1}}\C^{m\ast}
\rightarrow (S^N\C^{m\ast})^\ast = S^N\C^m.
$$
We have
$$
R(f_1, \dots , f_{m-1}, h) =
\langle Q(f_1, \dots , f_{m-1}), h^N \rangle
$$
for $(f_1, \dots , f_{m-1}, h)
\in S^{n_1}\C^{m\ast} \times \dots \times  S^{n_{m-1}}\C^{m\ast}
\times \C^{m\ast}$. It follows from this formula that
\roster
\item if $(f_1, \dots , f_{m-1}) \in
S^{n_1}\C^{m\ast} \times \dots \times  S^{n_{m-1}}\C^{m\ast}$,
then $Q(f_1, \dots , f_{m-1}) = 0$ iff
$dim(V(f_1) \cap \dots \cap V(f_{m-1})) \ge 1$,
\item if $(f_1, \dots , f_{m-1}) \in
S^{n_1}\C^{m\ast} \times \dots \times  S^{n_{m-1}}\C^{m\ast}$,
$|V(f_1) \cap \dots \cap V(f_{m-1})| = N$, then
$Q(f_1, \dots , f_{m-1}) = l_1 \cdot \dots \cdot l_N$,
there $l_i \in \C^m$ and
$$
V(f_1) \cap \dots \cap V(f_{m-1}) =
\{\overline {l}_1, \dots, \overline {l}_N \}.
$$
\endroster
It follows that if $(f_1, \dots , f_{m-1}) \in
S^{n_1}\C^{m\ast} \times \dots \times  S^{n_{m-1}}\C^{m\ast}$ and
$Q(f_1, \dots , f_{m-1}) \neq 0$, then
$Q(f_1, \dots , f_{m-1}) = l_1 \cdot \dots \cdot l_N$
and $\overline {l_i} \in V(f_1) \cap \dots \cap V(f_{m-1})$,
$1 \le i \le N$.
\par
Consider the group $SL_m \times SL_{m-1}$. The group
$SL_m \times SL_{m-1}$ acts canonically in the spaces
$\C^m, \C^{m-1}$, $\C^{m\ast}$, $S^n\C^{m\ast}$,
$S^n\C^{m\ast} \otimes \C^{m-1}$, \dots . Let
$a_1, \dots, a_{m-1}$ be the standart basis of $\C^{m-1}$.
\par
Consider the homogeneous (of degree $(m-1) n^{m-2}$) polynomial
mapping
$$
S^n\C^{m\ast} \otimes \C^{m-1} \rightarrow
S^{n^{m-1}}\C^m,
$$
$$
f_1 \otimes a_1 + \dots + f_{m-1} \otimes a_{m-1}
\mapsto Q(f_1, \dots , f_{m-1}).
$$
It follows from the previous notes that this mapping
is  $SL_m \times SL_{m-1}$-covariant.
\bigpagebreak
\par
\bf \S3. \rm
In this section we recall some facts about
representations of the groups $SL_2$ and $SL_3$.
\par
The group $SL_2$ acts canonically in the space
$\C^2$. Let $a_1, a_2$ be the standard basis
of $\C^2$ and $y_1, y_2$ be the dual basis of
$\C^{2\ast}$. The group $SL_3$ acts canonically in the space
$\C^3$. Let $e_1, e_2, e_3$ be the standard basis
of $\C^3$ and $x_1, x_2, x_3$ be the dual basis of
$\C^{3\ast}$.
\par
Recall that $SL_2$-module $V(b) = S^b\C^{2\ast}$
is irreducible and $dimV(b) = b+1$.
Every irreducible $SL_2$-module is isomorphic
to $V(b)$ for some $b$. For example,
$$
V(b) = S^b\C^{c\ast} \rightarrow S^b\C^2,
$$
$$
f(y_1, y_2) \mapsto f(a_2, -a_1)
$$
is the isomorphism of $SL_2$-moduli.
The decomposition
$$
V(b) \otimes V(b') \simeq V(b + b') \oplus V(b + b' - 2)
\oplus \dots \oplus V(b + b' - 2min\{b,b'\}).
$$
holds.
Therefore, there exists a unique (up to a scalar factor)
nontrivial bilinear covariant (transvectant)
$$
\psi_i : V(b) \times V(b') \rightarrow V(b + b' - 2i), \quad
0 \le i \le min\{b,b'\}.
$$
The explicite form of $\psi_i$ is
$$
\psi_i (f_1, f_2) =
\frac{(b-i)!(b'-i)!}{b!b'!}
\sum_{0 \le j \le i}
(-1)^j \binom ij
\frac{\partial^i f_1}{\partial y_1^{i-j} \partial y_2^j}
\frac{\partial^i f_2}{\partial y_1^j \partial y_2^{i-j}}.
$$
where $f_1 \in V(b)$, $f_2 \in V(b')$.
\par
Consider the linear $SL_3$-mapping
$$
\Delta = \sum \frac{\partial}{\partial e_i} \otimes
\frac{\partial}{\partial x_i} :
\C[\C^{3\ast}] \otimes \C[\C^3] \rightarrow
\C[\C^{3\ast}] \otimes \C[\C^3].
$$
For $b, c \ge 0$ set
$$
V(b,c) = Ker(\left.\Delta\right|_{S^b\C^3 \otimes S^c\C^{3\ast}}).
$$
Recall that $SL_3$-module $V(b,c)$ is irreducible and
$dimV(b,c) = \frac{1}{2}(b+1)(c+1)$ $(b+c+2)$.
Every irreducible $SL_3$-module is isomorphic to
$V(b, c)$ for some $b, c$. Let
$$
\pi_{b,c} : S^b\C^3 \otimes S^c\C^{3\ast}
\rightarrow V(b,c)
$$
be $SL_3$-projection.
The decomposition
$$
V(b,c) \otimes V(b',c') \simeq
V(b+b',c+c') \oplus V(b+b'+1,c+c'-2) \oplus
\dots
$$
$$
\oplus V(b+b'+min\{c,c'\},c+c'-2min\{c,c'\})
\oplus \dots .
$$
holds.
Therefore, there exists a unique (up to a scalar
factor) nontrivial bilinear covariant
$$
\rho_i : V(b,c) \times V(b',c')
\rightarrow V(b+b'+i,c+c'-2i), \quad
0 \le i \le min\{c,c'\}.
$$
Let us give the explicit form of $\rho_1$:
$$
\rho_1 (f_1, f_2) = \pi_{b+b'+1,c+c'-2}
\left(\sum_{\sigma \in S_3}
sgn(\sigma) e_{\sigma(1)}
\frac{\partial f_1}{\partial x_{\sigma(2)}}
\frac{\partial f_2}{\partial x_{\sigma(3)}}\right),
$$
where $f_1 \in V(b,c)$, $f_2 \in V(b',c')$.
Similarly, there exists a unique (up to a scalar
factor) nontrivial bilinear covariant
$$
\tau_i : V(b,c) \times V(b',c')
\rightarrow V(b+b'-2i,c+c'+i), \quad
0 \le i \le min\{b,b'\}.
$$
\proclaim{Lemma 3.1} Consider the polyhomogeneous
(of polydegree $(2n,1,1)$) covariant
$$
\eta : \C^{3\ast} \times S^n\C^{3\ast}
\times S^n\C^{3\ast}
\rightarrow V(2n-4,2),
$$
$$
(h,f_1,f_2) \mapsto \tau_2(\rho_n(h^n,f_1),\rho_n(h^n,f_2)).
$$
Then
$$
\eta (x_3,f_1,f_2) =
c_0 \left(\left.\psi_2(f_1(y_1,y_2,0),
f_2(y_1,y_2,0))\right|_{y_1 = e_2, y_2 = -e_1}\right) x_3^2,
$$
where $c_0 \in \C$, $c_0 \neq 0$.
\endproclaim
\demo{Proof}
Since $\rho_n \neq o$, we have
$$
\rho_n (x_3^n,f_i) = c' f_i(e_2,-e_1,0),
\ \ i = 1,2, \tag{3.1}
$$
where $c' \in \C$, $c' \neq 0$.
Since $\tau_2 \neq 0$, we have
$$
\tau_2(r_1(e_1,e_2,0),r_2(e_1,e_2,0)) =
$$
$$
c''\left(\left.\psi_2(r_1(-y_2,y_1,0),r_2(-y_2,y_1,0))
\right|_{y_1 = e_2, y_2 = -e_1}\right) x_3^2, \tag{3.2}
$$
where $r_i \in S^n\C^3$, $c'' \in \C$, $c'' \neq 0$.
By (3.1) and (3.2), it follows the Lemma.
\enddemo
\bigpagebreak
\par
\bf \S4. \rm
In this section we define group $G$, some
$G$-moduli, and covariants. We use the notations
of \S3.
\par
Let $h \in \C^{3\ast}, h \neq 0$. Define the
2-dimensional affine space
$$
\A (h) = P\C^3\setminus
\{\overline{x} \in P\C^3 |
\langle h,x \rangle = 0 \}.
$$
Define 2-form
$$
\omega_h = \left.\frac{\Omega_h}{h^2}\right|_{\A (h)}
$$
on $\A (h)$,
where $\Omega_h \in \wedge^2\C^{3\ast}$,
$\Omega_h \wedge dh =
dx_1 \wedge dx_2 \wedge dx_3.$
\par
Consider an affine space $\A (h)$, where
$h \in \C^{3\ast}$, $h \neq 0$.
Let $\C[\A (h)]_{\le n}$ be the linear space of
polynomials (of degree $\le n$) on $\A (h)$.
The space $\C[\A (h)]_{\le n} \otimes \C^2$
is the linear space of polynomial mappings
(of polynomial degree $\le n$) of affine
space $\A (h)$ to $\C^2$. Fix the following
isomorphisms of the linear spaces:
$$
i_h : S^n\C^{3\ast} \rightarrow \C[\A (h)]_{\le n},
$$
$$
f(x) \mapsto \frac{f(x)}{(h(x))^n},
$$
$$
I_h : S^n\C^{3\ast} \otimes \C^2
\rightarrow \C[\A (h)]_{\le n} \otimes \C^2,
$$
$$
f_1(x) \otimes a_1 + f_2(x) \otimes a_2
\mapsto \frac{f_1(x)}{(h(x))^n} \otimes a_1 +
\frac{f_2(x)}{(h(x))^n} \otimes a_2.
$$
\par
Set
$$
G = SL_3 \times SL_2.
$$
The group $G$ acts canonically in the spaces
$\C^3, \C^2, \C^{3\ast}$, $S^n\C^{3\ast} \otimes \C^2, \dots$.
Consider the representation
$$
G: \C^{3\ast} \times S^n\C^{3\ast} \otimes \C^2.
$$
\par
The decomposition
$$
\C^{3\ast} \otimes S^2(S^n\C^{3\ast} \otimes \C^2))
\simeq S^{2n-2}\C^{3\ast} \oplus \dots.
$$
holds.
Therefore, there exists a unique (up to a scalar
factor) nontrivial polyhomogeneous
(of polydegree (1,2)) covariant
$$
J : \C^{3\ast} \times S^n\C^{3\ast} \otimes \C^2
\rightarrow S^{2n-2}\C^{3\ast}.
$$
Let us give the explicit form of $J$:
$$
J(h,f_1(x) \otimes a_1 + f_2(x) \otimes a_2)=det
\left(\matrix\frac{\partial f_1}{\partial x_1} &
             \frac{\partial f_1}{\partial x_2} &
             \frac{\partial f_1}{\partial x_3} \\
             \frac{\partial f_2}{\partial x_1} &
             \frac{\partial f_2}{\partial x_2} &
             \frac{\partial f_2}{\partial x_3} \\
             \frac{\partial h}{\partial x_1} &
             \frac{\partial h}{\partial x_2} &
             \frac{\partial h}{\partial x_3} \endmatrix \right).
$$
\proclaim{Lemma 4.1} Suppose
$(h,f) \in \C^{3\ast} \times S^n\C^{3\ast} \otimes \C^2$,
$h \neq 0$; then
$i_h(J(h,f))$ is the jacobian of the polynomial
mapping $I_h(f)$.
\endproclaim
\demo{Proof} It is sufficient to consider the case
$h(x) = x_3$. The form $\omega_h$ on $\A (h)$ is
$\frac{dx_1 \wedge dx_2}{x_3^2}$. We have
$$
I_h(f_1(x) \otimes a_1 + f_2(x) \otimes a_2)(\overline x) =
\frac{f_1(x)}{x_3^n}a_1 + \frac{f_2(x)}{x_3^n}a_2,
$$
$$
d(\frac{f_1(x)}{x_3^n}) \wedge d(\frac{f_2(x)}{x_3^n}) =
\frac{1}{x_3^{2n-2}}
det \left(\matrix\frac{\partial f_1}{\partial x_1} &
                 \frac{\partial f_1}{\partial x_2} \\
                 \frac{\partial f_2}{\partial x_1} &
                 \frac{\partial f_2}{\partial x_2} \endmatrix
    \right) \frac{dx_1 \wedge dx_2}{x_3^2},
$$
$$
i_h^{-1}(\frac{1}{x_3^{2n-2}}
det \left(\matrix\frac{\partial f_1}{\partial x_1} &
                 \frac{\partial f_1}{\partial x_2} \\
                 \frac{\partial f_2}{\partial x_1} &
                 \frac{\partial f_2}{\partial x_2} \endmatrix
    \right)) =
det \left(\matrix\frac{\partial f_1}{\partial x_1} &
                 \frac{\partial f_1}{\partial x_2} \\
                 \frac{\partial f_2}{\partial x_1} &
                 \frac{\partial f_2}{\partial x_2} \endmatrix
    \right) =
$$
$$
J(x_3,f_1(x) \otimes a_1 + f_2(x) \otimes a_2).
$$
\enddemo
Set
$$
J_c : \C^{3\ast} \times S^n\C^{3\ast} \otimes \C^2
\rightarrow V(1,2n-3),
$$
$$
(h,f)
\mapsto \rho_1(h,J(h,f)).
$$
$J_c$ is the covariant of polydegree (2,2).
\proclaim{Lemma 4.2} Suppose
$(h,f) \in \C^{3\ast} \times S^n\C^{3\ast} \otimes \C^2$,
$h \neq 0$; then
$J_c(h,f) = 0$ iff the jacobian of the polynomial
mapping $I_h(f)$  is equal to a constant.
\endproclaim
\demo{Proof} Suppose the jacobian of the polynomial
mapping $I_h(f)$ is equal to a constant, i.e.
$J(h,f) = ch^{2n-2}$, $c \in \C$ (Lemma 4.1).
We have
$$
J_c(h,f) = \rho_1(h,J(h,f)) = \rho_1(h,ch^{2n-2}) = 0.
$$
\par
Suppose $J_c(h,f) = 0$. From
$$
0 = J_c(h,f) = \rho_1(h,J(h,f))
$$
it follows that $J(h,f) = ch^{2n-2}$,
$c \in \C$ and thus the jacobian
of the polynomial mapping $I_h(f)$ is equal to a constant
(Lemma 4.1).
\enddemo
\bigpagebreak
\par
\bf \S5. \rm
In this section we formulate conjecture 1.
This conjecture is true iff the 2-dimensional
jacobian conjecture is
true. Conjecture 1 is a partial case of general
invariant theretical problem * (see \S1).
We use the notations of \S3 and \S4.
\par
Set
$$
D_i : \C^{3\ast} \times S^n\C^{3\ast} \otimes \C^2
\rightarrow S^{n^2-i}\C^3,\ \ \ i \ge 0,
$$
$$
(h,f_1(x) \otimes a_1 + f_2(x) \otimes a_2)
\mapsto \Delta^i(h^iQ(f_1,f_2)),
$$
(the definition of the covariant $Q$ see
in \S2). $D_i$ is the covariant of polydegree
$(i,n^2)$.
\par
Let $(h,f) = (h,f_1(x) \otimes a_1 + f_2(x) \otimes a_2)
\in \C^{3\ast} \times S^n\C^{3\ast} \otimes \C^2$,
$h \neq 0$. A polynomial mapping $I_h(f)$ is called
in general position if $Q(f_1,f_2) \neq 0$.
It is obvious that a polynomial mapping
$I_h(f)$ is not in general position iff
$dim(I_h(f)^{-1}(0)) \ge 1$ or
the polynomial degree of $I_h(f)$ is less
than $n$.
\proclaim{Lemma 5.1} Suppose
$(h,f) = (h,f_1(x) \otimes a_1 + f_2(x) \otimes a_2)
\in \C^{3\ast} \times S^n\C^{3\ast} \otimes \C^2$,
$h \neq 0$, and $i \ge 0$; then
$D_i(h,f) = 0$ iff the polynomial mapping
$I_h(f)$ is not in general position or
$|I_h(f)^{-1}(0)| < i$.
\endproclaim
\demo{Proof} It is obvious that the lemma
is true if the polynomial
mapping $I_h(f)$ is not in general position.
Suppose the polynomial mapping $I_h(f)$ is in
general position.
If $|I_h(f)^{-1}(0)| < i$, then
$$
Q(f_1,f_2) = l_1 \cdot \dots \cdot l_{i-1}
\cdot \overline{l}_1 \cdot \dots
\cdot \overline{l}_{n^2-i+1},
$$
where $\langle h, \overline{l}_1 \rangle = 0$,
$1 \le j \le n^2 - i + 1$. We have
$$
D_i(h,f) = \Delta^i(h^iQ(f_1,f_2))=
$$
$$
\Delta^i(h^i l_1 \cdot \dots \cdot l_{i-1}
\cdot \overline{l}_1 \cdot \dots
\cdot \overline{l}_{n^2-i+1}) = 0.
$$
\par
If
$$
0 = D_i(h,f) = \Delta^i(h^iQ(f_1,f_2)),
$$
then
$$
Q(f_1,f_2) = l_1 \cdot \dots \cdot l_{i-1}
\cdot \overline{l}_1 \cdot \dots
\cdot \overline{l}_{n^2-i+1},
$$
where $\langle h, \overline{l}_j \rangle = 0,
1 \le j \le n^2 - i + 1$ and thus
$|I_h(f)^{-1}(0)| < i$.
\enddemo
\proclaim{Conjecture 1} Consider the $G$-module
$\C^{3\ast} \times S^n\C^{3\ast} \otimes \C^2$
and the covariants $J_c, D_2$; then
$$
\left.D_2(h,f)\right|_{J_c(h,f)=0} = 0.
$$
\endproclaim
It follows from Lemma 4.2 and Lemma 5.1 that
conjecture 1 is true iff the following
conjecture is true.
\proclaim{Conjecture $1'$} Suppose the jacobian
of a polynomial mapping
$$
F : \C^2 \rightarrow \C^2
$$
is equal to a constant; then
$dimF^{-1}(0) \ge 1$ or
$|F^{-1}(0)| \le 1$.
\endproclaim
Evidently, conjecture $1'$ is true iff the 2-dimensional
Jacobian conjecture is true.
\bigpagebreak
\par
\bf \S6. \rm
In this section we formulate conjecture 2.
This conjecture is true iff the 2-dimensional
Jacobian conjecture is true.
Conjecture 2 is a partial case of general
invariant theretical problem ** (see \S1).
We use the notations of \S3 and \S4.
\proclaim{Conjecture 2}
Consider the linear representation
$G : \C^{3\ast} \times S^n\C^{3\ast} \otimes \C^2$
and the elements
$$
(x_3,f_0) \in \C^{3\ast} \times S^n\C^{3\ast} \otimes \C^2, \quad
x_3^{2n-2} \in S^{2n-2}\C^{3\ast},
$$
where
$$
f_0 = x_1^{n_1} x_2^{n_2} \otimes a_1 +
      x_1^{n_1} x_2^{n_2} \otimes a_2,
$$
$n_1 \ge 1$, $n_2 \ge 1$, $n_1 + n_2 = n$;
then there exists a polyhomogeneous
$(J,J_c)$-identity
$$
\gamma : \C^{3\ast} \times S^n\C^{3\ast} \otimes \C^2
\times S^{2n-2}\C^{3\ast} \rightarrow W
$$
such that
$$
\gamma (x_3,f_0,x_3^{2n-2}) \neq 0.
$$
\endproclaim
\proclaim{Proposition 6.1}
Conjecture 2 is true iff the 2-dimensional
Jacobian conjecture is true.
\endproclaim
\demo{Proof} \bf1. \rm Let us deduce the 2-dimensional
Jacobian conjecture from conjecture 2.
\par
If the 2-dimensional Jacobian conjecture is not true,
then there exists a counterexample
$$
F = (F_1,F_2) : \C^2 \rightarrow \C^2,
$$
$$
z = (z_1,z_2) \mapsto (F_1(z),F_2(z))
$$
such that
$$
F_1(z) = z_1^{n_1}z_2^{n_2} + \widetilde {F}_1(z), \quad
F_2(z) = z_1^{n_1}z_2^{n_2} + \widetilde {F}_2(z),
$$
where $n_1 \ge 1$, $n_2 \ge 1$, $n_1 + n_2 = n$,
$deg \widetilde {F}_1 < n$,
$deg \widetilde {F}_2 < n$ (Folk-theorem).
Consider the linear representation
$G : \C^{3\ast} \times S^n\C^{3\ast} \otimes \C^2$
and the elements
$$
(x_3,f_0) \in \C^{3\ast} \times S^n\C^{3\ast} \otimes \C^2, \quad
x_3^{2n-2} \in S^{2n-2}\C^{3\ast},
$$
where
$$
f_0 = x_1^{n_1} x_2^{n_2} \otimes a_1 +
      x_1^{n_1} x_2^{n_2} \otimes a_2.
$$
It follows from conjecture 2 that there exists
a polyhomogeneous $(J,J_c)$-identity
$$
 \gamma : \C^{3\ast} \times S^n\C^{3\ast} \otimes \C^2
\times S^{2n-2}\C^{3\ast} \rightarrow W
$$
such that
$$
\gamma (x_3, f_0, x_3^{2n-2}) \neq 0.
$$
\par
Set
$$
f = x_3^n F_1(\frac{x_1}{x_3},\frac{x_2}{x_3}) \otimes a_1 +
    x_3^n F_2(\frac{x_1}{x_3},\frac{x_2}{x_3}) \otimes a_2 =
f_0 + \widetilde {f},
$$
where
$$
\widetilde {f} = x_3^n \widetilde {F}_1(\frac{x_1}{x_3},\frac{x_2}{x_3})
\otimes a_1 +
x_3^n \widetilde {F}_2(\frac{x_1}{x_3},\frac{x_2}{x_3})
\otimes a_2.
$$
By Lemma 4.1, it follows that
$J(x_3,f) = x_3^{2n-2}$.
Therefore, $J_c(x_3,f) = 0$ and
$$
0 = \gamma (x_3,f,J(x_3,f)) =
\gamma (x_3, f_0 + \widetilde {f}, x_3^{2n-2}). \tag{6.1}
$$
Let $(d_1,d_2,d_3)$ be the polydegree of the covariant
$\gamma$. Set
$$
N = 2 d_1 - n d_2 + 2 (2n-2) d_3,
$$
$$
g(t) = (\left(\matrix t^{-1} &  0      & 0   \\
                       0     & t^{-1}  & 0   \\
                       0     &  0      & t^2 \endmatrix
        \right),
        \left(\matrix 1 & 0 \\
                      0 & 1 \endmatrix
        \right)) \in G.
$$
{}From (6.1) we obtain
$$
0 = t^N g(t) \cdot
\gamma (x_3, f_0 + \tilde {f}, x_3^{2n-2}) =
$$
$$
\gamma (t^2 g(t) \cdot x_3,
t^{-n} (g(t) \cdot f_0 + g(t) \cdot \tilde {f}),
t^{2(2n-2)} g(t) \cdot x_3^{2n-2})
$$
$$
\rightarrow \gamma (x_3, f_0, x_3^{2n-2}) \neq 0
$$
as $t \rightarrow \infty$. This contradiction
concludes the proof.
\par
\bf 2. \rm Let us deduce conjecture 2 from
the 2-dimensional Jacobian conjecture.
\par
Set
$$
\widetilde \eta : \C^{3\ast} \times S^n\C^{3\ast} \otimes \C^2
\rightarrow V(2n-4,2) \otimes S^2\C^2,
$$
$$
(h, f_1 \otimes a_1 + f_2 \otimes a_2)
\mapsto \eta (h,f_1,f_1) \otimes a_1^2 +
      2 \eta (h,f_1,f_2) \otimes a_1 a_2 +
        \eta (h,f_2,f_2) \otimes a_2^2
$$
(the definition of the covariant $\eta$
see in \S3). $\widetilde {\eta}$ is
the polyhomogeneous (of polydegree
$(2n,2)$) covariant.
\par
Suppose $J(x_3,f) = cx_3^{2n-2}$,
where $c \in \C$, $c \neq 0$,
$f = f_1 \otimes a_1 + f_2 \otimes a_2$; then
by Lemma 4.1 and 2-dimensional \ Jacobian
\ cojecture, \ it \ follows \ that \
$f_1(x_1,x_2,0) = c_1 \lambda (x_1,x_2)^n$,
$f_2(x_1,x_2,0) = c_2 \lambda (x_1,x_2)^n$,
where
$\lambda (x_1,x_2)$ is a linear form in the variables
$x_1,x_2$, and by Lemma 3.1, it follows that
$\widetilde {\eta} (x_3,f) = 0$.
Therefore,
$$
\left.\widetilde {\eta}
(x_3,f)\right|_{J(x_3,f) - c x_3^{2n-2} = 0} = 0 \tag{6.2}
$$
for $c \in \C$, $c \neq 0$.
\par
By (6.2), it follows that
$$
\left.\widetilde {\eta}
(h,f)\right|_{J(h,f) - c h^{2n-2} = 0} = 0 \tag{6.3}
$$
for $c \in \C$, $c \neq 0$.
\par
Consider the homogeneous (of polydegree (2n,2,1))
covariant
$$
\gamma : \C^{3\ast} \times S^n\C^{3\ast} \otimes \C^2
\times S^{2n-2} \C^{3\ast}
\rightarrow (V(2n-4,2) \otimes S^2 \C^2)
\otimes S^{2n-2} \C^{3\ast}
$$
$$
(h,f,j) \mapsto \widetilde {\eta} (h,f) \otimes j.
$$
By Lemma 4.2 and (6.3) it follows that
$$
\left.\gamma(h,f,J(h,f))\right|_{J_c(h,f) = 0} = 0.
$$
This means that $\gamma$ is a polyhomogeneous
$(J,J_c)$-identity.
Using Lemme 3.1, we get
$$
\widetilde \eta (x_3,f_0) =
\eta (x_3,x_1^{n_1}x_2^{n_2},x_1^{n_1}x_2^{n_2}) \otimes
(a_1^2 + 2 a_1 a_2 + a_2^2) =
$$
$$
c_0 \frac{n_1n_2(1-n_1-n_2)}{n^2(n-1)^2}
e_2^{2n_1 - 2} e_1^{2n_2 - 2} x_3^2 \otimes
(a_1^2 + 2 a_1 a_2 + a_2^2) \neq 0.
$$
Therefore,
$$
\gamma (x_3, f_0, x_3^{2n-2}) =
\widetilde \eta (x_3,f_0) \otimes x_3^{2n-2} \neq 0.
$$
\enddemo


\Refs
\ref \no1
\by W. Fulton and J.Harris
\book Representation Theory, {\rm
Springer Graduate Texts in Math., v. 129}
\publ Springer-Verlag, Berlin
\yr 1991
\endref
\ref \no2
\by H. Bass, E.H. Connell, and D. Wright
\paper The Jacobian conjecture:
reduction of degree and formal expansion of
the inverse.
\jour Bull. Amer. Math. Soc.
\vol 1982 \issue 7 \pages 287--330
\endref
\ref \no3
\by O. Keller
\paper Ganze Cremona-Transformationen
\jour Monatsh. Math. Phys.
\vol 1939 \issue 47 \pages 299--306
\endref
\endRefs
\enddocument